\documentclass[journal]{IEEEtran}
\usepackage{multirow}
\usepackage{multicol}
\usepackage{lipsum}
\usepackage{graphicx}
\usepackage{graphics}
\usepackage{amsmath}
\usepackage{amsbsy}
\usepackage{blindtext}
\usepackage{tcolorbox}
\usepackage{amssymb}
\usepackage{scalerel}
\usepackage{mathabx}
\usepackage{verbatim}
\usepackage{booktabs}
\usepackage{rotating,tabularx}
\usepackage{mathrsfs} 
\usepackage{hyperref}

\usepackage{adjustbox}
\usepackage{soul}
\usepackage{xcolor}
\usepackage{cite}
\usepackage{tikz}
\usepackage{color, colortbl}
\usepackage[first=0,last=9]{lcg}
\definecolor{LightGray}{rgb}{0.7,0.7,0.7}

\usepackage[wby]{callouts}
\usetikzlibrary{shapes.multipart}

\usepackage{array}
\usepackage{makecell}

\allowdisplaybreaks

\usepackage{amsthm}
\theoremstyle{definition}

\theoremstyle{remark}

\usepackage[utf8]{inputenc}
\usepackage[english]{babel}

\usepackage[caption=false,font=footnotesize]{subfig}
\usepackage[T1]{fontenc}
\usepackage{scalerel,stackengine}
\newcommand\reallywidecheck[1]{%
\savestack{\tmpbox}{\stretchto{%
  \scaleto{%
    \scalerel*[\widthof{\ensuremath{#1}}]{\kern-.6pt\bigwedge\kern-.6pt}%
    {\rule[-\textheight/2]{1ex}{\textheight}}
  }{\textheight}%
}{0.5ex}}%
\stackon[1pt]{#1}{\scalebox{-1}{\tmpbox}}%
}
\stackMath
\IEEEoverridecommandlockouts
\IEEEoverridecommandlockouts

\renewcommand{\baselinestretch}{0.991}

\newcommand*{\rn}{\textcolor{black}}

\newif\ifarxiv
\arxivtrue
\arxivfalse

\begin{document}

\title{\LARGE\bf
\rn{AC False Data Injection Attacks in Power Systems: Design and Optimization}}

\author{Mohammadreza Iranpour$^{\ast}$, Mohammad Rasoul Narimani$^{\ast}$
\thanks{${\ast}$: Department of Electrical and Computer Engineering, California State University Northridge (CSUN). Rasoul.narimani@csun.edu. Support from NSF contract \#2308498.}%
}

\maketitle

\begin{abstract}
\rn{False Data Injection (FDI) attacks are one of the challenges that the modern power system, as a cyber-physical system, is encountering. Designing AC FDI attacks that accurately address the physics of the power systems could jeopardize the security of power systems as they can easily bypass the traditional Bad Data Detection (BDD) algorithm. Knowing the essence of the AC FDI attack and how they can be designed gives insight about detecting the system again these attacks. Moreover, recognition of the nature of these attacks, especially when they are designed optimally, is essential for benchmarking various defensive approaches to increase the resilience of power systems. This paper presents a unified approach to demonstrate the process of designing optimal AC FDI attack. In this connection, we first define the process of designing an AC-based FDI attack that satisfies AC power flow equations. We then formulate an optimization problem to design an optimal AC FDI attack that both satisfies AC power flow equations and overloads a specific line in the system. The objective function is defined to optimize the magnitude of the attack vector in such a way that it can evade residue-based BDD approaches. The proposed approach for designing AC FDI attacks is applied to the IEEE 118-bus test case system. Various comparisons are conducted to elaborate on the impact of optimally designing AC FDI attacks on the residual for the AC state estimation algorithm. Comparing the results of optimal and non-optimal AC FDI attacks demonstrates the impact on the difficulty of detecting FDI attacks and the importance of optimally designing these attacks.} 
    

\end{abstract}

\section{Introduction}
\label{Introduction}
\rn{Although cyber-physical networks have the potential to enhance the efficiency, controllability, and reliability of power systems by leveraging the advantages of Information and Communication Technologies (ICT), these systems also introduce vulnerabilities. Specifically, the integration of communication systems into power systems can jeopardize their security, as attackers may hack measurement and communication devices to intrude false data into power system measurements, raising concerns about the secure performance of these networks~\cite{rahman2013false}. False Data Injection (FDI) attacks, which involve the deliberate insertion of incorrect data into a system's measurements, have emerged as a significant threat to the security of power grids. These attacks pose a major challenge to the use of ICT applications in power systems~\cite{rahman2013false, boyaci2022infinite, boyaci2022generating}.}

\rn{The state estimator, a core component of the Supervisory Control and Data Acquisition (SCADA) system, continuously monitors the operating state of a power system to ensure safe and reliable operation. However, it has been demonstrated that intruders can hack multiple sensors and measurement devices, compromising the state estimator by injecting pre-designed false data vectors into the meters which is known as the FDI attacks~\cite{liu2011false, boyaci2021graph, boyaci2021joint}. If this attack vector meets certain conditions, it can bypass the commonly used residue-based Bad Data Detectors (BDD), falsifying network topology and misleading the control center. This undetected manipulation can lead to various disruptive consequences, such as intentional branch outages or power grid frequency excursions, which may result in blackouts or damage to electrical equipment~\cite{du2021targeted, boyaci2022spatio}.}  





\rn{Over the past decade, numerous studies have reviewed FDI attacks, covering background materials, construction methods, and detection and defense strategies \cite{zhang2019false}. The construction methods have been examined from various perspectives, including FDI attacks with limited budgets, those based on state estimation with incomplete system knowledge\cite{zhao2018generalized}, and those using a data-driven approach with incomplete system knowledge\cite{zhang2018can}. The focus of FDI attacks with limited budgets is on constructing valid attacks with minimal effort, such as compromising at most k meters, known as the k-sparse problem\cite{liu2011false}. In this scenario, the attack designer needs to formulate an optimal strategy tailored to specific situations and targets.}


\rn{Building on this background, two realistic categories of FDI attacks have been considered from a design perspective. The first is random false data injection attacks, where the attacker aims to find any attack vector that can lead to incorrect estimation of state variables. The second is targeted false data injection attacks, where the attacker aims to find an attack vector that can inject a specific error into certain state variables. Research by Liu et al. \cite{liu2011false} demonstrates that attackers can systematically and efficiently construct attack vectors for both scenarios. Additionally, in \cite{jin2018power} is shown that these systematic approaches can be implemented by formulating an FDI attack as an optimization problem, to find a stealthy and sparse data injection vector on the sensor measurements to produce spurious and misleading state estimates.}

\rn{Numerous studies have focused on optimization approaches that emphasize sparsity and minimizing the number of compromised measurements, aiming to incorporate constraints related to limited resources and stealthy operations \cite{jin2017semidefinite,liu2020network}. Other studies have analyzed these optimization approaches from different perspectives. For instance, semi-definite programming (SDP) relaxation and sparsity penalties have been used to recover a near-optimal solution \cite{jin2017semidefinite}. A reduced row-echelon (RRE) form-based greedy method has been employed to compute the minimum number of targeted measurements needed for an attack \cite{nayak2020modelling}. Additionally, genetic algorithms and neural networks have been applied to construct the least-effort attack vector \cite{jin2018power}. The minimum number of sensors required to stage stealth FDI attacks has also been quantified by formulating a minimum cardinality problem, with various algorithms proposed for efficient computation \cite{jin2018power}.}

\rn{In addition to the various approaches to designing an optimal False Data Injection (FDI) attack, another perspective involves optimizing the magnitudes of the attack vector. In this scenario, attackers aim to adjust the magnitudes of injected false data to ensure they fall within acceptable bounds, evading detection by traditional state estimation algorithms.} \rn{This paper addresses this aspect by formulating an objective function to optimize the magnitudes of injected false data in power flow measurements within the power system. Specifically, it focuses on AC-based FDI attacks, known for their complexity in design, and difficulty in detection compared to DC-based FDI attacks \cite{rahman2013false,tran2021designing,wang2013cyber}. Using the nonlinear AC power flow model to accurately represent power flow physics results in a non-convex AC power flow, increasing the complexity of the problem~\cite{narimani2023tightening,narimani2020tightening, narimani2020strengthening,narimani2018comparison, narimani2018empirical,narimani2018improving}. These attacks systematically leverage non-linear power flow equations and Kirchhoff's laws in a compromised power system.}
\rn{Apart from defining an objective function for optimizing attack vector magnitudes, this paper also considers the non-linear power flow equations of the compromised system as constraints within the. This simultaneous optimization of a successful FDI attack while optimizing the attack vector magnitudes enhances the stealthiness of the attack, making it challenging to detect using existing methods.}

\rn{The remainder of this paper is organized as follows: In Section~\ref{sec:Preliminary}, we present a systematic approach to designing a successful AC-based False Data Injection (FDI) attack, outlining the definition of attacked and non-attacked zones within the power system and corresponding power flow equations. Section~\ref{sec:Optimal FDIA} formulates an appropriate objective function and necessary constraints for designing an optimal AC-based FDI attack. The implementation and analysis of this approach on the IEEE 118-bus test system are discussed in Section~\ref{sec:ANALYSIS}. Finally, Section~\ref{sec:conclusion} offers concluding remarks on the study's findings and implications.}


\section{Preliminary Concepts for AC False Data Injection Attack in Power Systems}
\label{sec:Preliminary}

\rn{To ensure continuous operation of power systems, engineers use monitoring systems that provide measurement data to state estimation platforms in control centers. State estimators determine various parameters, such as voltage levels and angles at buses, enabling appropriate decision-making and actions. Fundamentally, accurate sensor readings usually yield state variable estimations that closely match their true values, while abnormal measurements can greatly alter these estimates. To detect these inconsistencies, researchers in power systems suggest calculating the measurement residual $r= z-H(x)$, which demonstrates the difference between the observed measurements ($z$) and estimated measurements ($h(x)$). When $\|z - H(x)\|>\tau$, where $\tau$ is a specified threshold, it is concluded that there are erroneous measurements.}

\rn{Although Bad Data Detectors (BDD) are engineered to detect and eliminate measurements that appear incorrect or fall outside the normal data distribution range, adversaries can inject meticulously crafted measurements that disrupt the data distribution while eluding detection. Assuming the residue of state-estimated variables under attack is: }

\begin{align}
&\nonumber  r_{attack}=\| z_{attack} - h(x_{attack}) \|  \\
&\nonumber = \| z_{attack}  - h(x_{attack}) + h(x) - h(x) \| \\
&\nonumber = \| (z +a - h(x_{attack}) + h(x) - h(x)  \| \\
& =  \| (r +a  - h(x_{attack}) + h(x)  \|.
\label{eq:attack vector}
\end{align}
    
\rn{\noindent where the attack vector $a$ is defined as:}

\begin{align}
  a = \| h(x_{attack}) - h(x) \|.
\label{eq:vector a}
\end{align}
 
\rn{The residue-based bad
data detection tests cannot detect the attack vector $a$ since the injected false data no longer affects the residue.}

\rn{The primary obstacle in developing AC attacks involves defining the function $h(x_{attack})$. This task requires certain foundational assumptions which, when considered, allow for defining this function and executing a successful attack. $h$ is the set of nonlinear power flow equations that relate measurement and variable states. To design a successful FDI attack, the non-linear relationship in power flow equations that connects state variables with various quantities, including power flow in lines and power injections into buses, should be considered. In the following sections the process of determining function $h(x_{attack})$, is presented.}

\subsection{Identifying Attack Zone}
\label{subsec:attack zone}

\rn{In general, an AC FDI attack focuses on a particular sub-grid region, leading to the division of the entire grid into two distinct areas: the ``attack area'' and the ``normal area''. The attack area encompasses the region directly targeted by the attack and consists of a set of buses forming a closed region where all alterations due to the attack occur internally. Conversely, the normal area remains unaffected by the attack.} 

\rn{As a fundamental aspect in designing a successful FDI attack to circumvent the BDD algorithm, every modification resulting from such an attack must be justified.
The initial assumption crucial for designing such an attack is to meet specific criteria between the normal and attacked areas, considering the multiple power transfers occurring between the affected and unaffected regions. Consequently, during the execution of an FDI attack on the targeted region, it becomes imperative to maintain equivalent total power transfers between this area and the unaffected regions. According to \cite{hug2012vulnerability}, fulfilling this requirement entails enclosing the attack area with buses equipped for power injection, i.e. buses with non-zero power injection, thus justifying every power change solely within the attack zone. The ability to justify power changes is attributed to buses equipped with power injection measurement devices, because they can justify the power changes through the crafting of injection measurements, as opposed to zero injection buses, which only extend the attack zones's influence.  Furthermore, to confine all alternations within the attack area during an FDI attack, it is essential to keep the state variables on the boundary unchanged.} 

\begin{figure}
    \centering
\includegraphics[scale=0.45,trim= 6cm 7.1cm 6.0cm 3.0cm,clip]{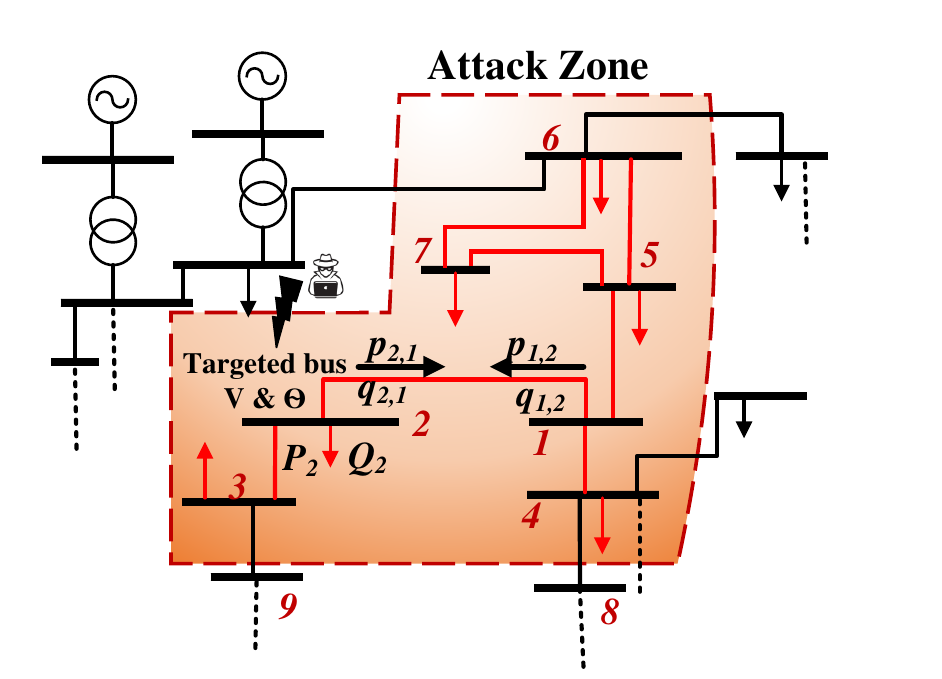}
	\caption{\rn{A toy example of defining the attack zone in a sample test case. The red line indicates the branches within the attack zone. The red dashed line separates the attack zone from the rest of the system.}}
	\label{fig:sample}
\end{figure}

\rn{In general, determining the attack zone involves designating specific buses, either individually or as groups, as primary focal points. These focal buses are then expanded by considering neighboring non-zero injection busses and passing through zero-injection buses to delineate the boundary of the attack zone. 
For example, consider Figure~\ref{fig:sample}, where bus $2$ is the targeted bus for attack. Altering the state variables at bus $2$ affects the power flow between busses $2$ and $1$. Since bus $1$ is a zero-injection bus, these alternations propagate to other links connected to it, expanding the attack effects. Similarly, the alternation of power flow between buses $2$ and $3$ occurs, with bus $3$ being a non-zero injection bus. Bus $3$ is a no-zero injection bus and this alternation of power flow can be justified by power injection measurement of this bus. By fixing the state variables (voltage magnitude and angle) at bus $3$, the power flow does not alter on the other link that is connected to this bus, thereby confining the attack region and preventing its expansion to other buses.  This process is repeated for buses $4$ and $6$, where the attack effects are bonded. In Figure~\ref{fig:sample}, the power flow and power injection measurements that need alternation due to the expanding attack effects are highlighted in red.}  \rn{In the following sections, the second assumption is presented, outlining a unified method for generating constraints to provide function $h(x)$ and subsequently designing and executing an AC FDI attack.} 


\subsection{Constraints-based AC False Data Injection Attack}
\label{subsec: AC FDIA}

\rn{By considering this assumption that the algebraic sum of the generated and consumed power in the attack region must remain unchanged, it is possible to establish two types of constraints that can define function $h(x)$, for designing the attack vector. In this regard:}

\begin{enumerate}
    \item \rn{In zero-injection buses within the attack zone, the algebraic sum of active/reactive power flows after the attack must be equal to zero. For example in the Figure~\ref{fig:sample}, since bus $1$ is a zero-injection bus, we can express this as: $\tilde{p}_{1,2}+ \tilde{p}_{2,1} + \tilde{p}_{1,5} + \tilde{p}_{5,1} + \tilde{p}_{1,4} +\tilde{p}_{4,1} = 0$ and $\tilde{q}_{1,2}+ \tilde{q}_{2,1} + \tilde{q}_{1,5} + \tilde{q}_{5,1} + \tilde{q}_{1,4} +\tilde{q}_{4,1} = 0$. In these equations, $\tilde{p}_{i,j}$ and $\tilde{q}_{i,j}$ represent the active and reactive power flow from bus $i$ to bus $j$, while $\tilde{p}_{j,i}$ and $\tilde{q}_{j,i}$ represent the active and reactive power flow for to bus $i$ from bus $j$. }

    \item \rn{In non-zero injection busses within the attack zone, the power injection after the attack equals the primary injection power plus the sum of all changes in power flows of lines connected to this non-zero injection bus and existing in the attack zone. For example in Figure~\ref{fig:sample}, since bus $4$ is a nonzero-injection bus, we can express its power flow injection after the attack as follows:$\tilde{P}_{4} = {P}_{4} + (\tilde{p}_{1,4} + \tilde{p}_{4,1} - {p}_{1,4} - {p}_{4,1})$ and $\tilde{Q}_{4} = {Q}_{4} + (\tilde{q}_{1,4} + \tilde{q}_{4,1} - {q}_{1,4} - {q}_{4,1})$. In these equations, $P_{4}$ and $Q_{4}$ represent the active and reactive power injections at the bus $4$. It is important to note that in equations concerning the boundary non-zero injection busses, only changes in the power flow of lines within the attack zone and connected to these busses should be considered. Eventually, based on these two assumptions and by solving these equations, the attack vector $a$, as per Equation~\eqref{eq:attack vector}, can be calculated.}
\end{enumerate}



\section{Optimal Designing of False Data Injection Attack}
\label{sec:Optimal FDIA}

\rn{To increase the chance of the FDI attack in bypassing the BDD algorithms, especially the residue-based BDD approaches, it is important to design the attacked vector properly. 
The main contribution of this paper is to optimize the elements of the attack vectors in order to have less impact on the residuals. By optimally designing the AC FDI attack, we not only design an AC FDI attack but also optimize its magnitude to make it more stealthy. In this regard, and based on the descriptions in the section~\ref{sec:Preliminary}, We can design the attack vector $a$ in Equation~\eqref{eq:vector a} by optimizing its elements. That is, we inject values into the measurement to minimize its impact on the AC state estimator. This minimization entails reducing the difference between corrupted measurements and the real measurements, as indicated in Equation~\eqref{eq:optimization}.}

\rn{Suppose the sets of buses, and lines of a system are demonstrated by $\mathcal{B}$, and $\mathcal{L}$, respectively. Also, $\mathcal{B_A}$, and $\mathcal{L_A}$, are corresponding sets of buses and lines in the attack zone,  respectively. Let $S_m = {P}_m + j {Q}_m$ represents the complex power injection, $V_m$ and $\theta_m$ represent the voltage magnitude and angle at bus~$m\in\mathcal{B_A}$, each line $\left(m,l\right)\in\mathcal{L_A}$ is modeled as a $\Pi$ circuit with mutual admittance $g_{ml}+j b_{ml}$ and shunt admittance $j b_{c,ml}$ and the voltage angle difference between buses $m$ and $l$ for $(m,l)\in\mathcal{L}$ is denoted as $\theta_{ml}=\theta_{m}-\theta_{l}$.
The difference between state variables, including both voltage magnitudes and angles and their values prior to conducting FDI attack is shown in Equation~\eqref{eq:optimization} as vector $c$.}

\begin{equation}
\label{eq:optimization}
 c= [\tilde{V}_{m}- V_{m,fix}, \tilde{\theta}_{m}-\theta_{m,fix}]
\end{equation}

\rn{Where $\tilde{V}_{m}$ and $\tilde{\theta}_{m}$ are variables representing voltage magnitude and angle that need to be calculated to conduct the attack, and $V_{m,fix}$ and $\theta_{m,fix}$ are the known values of voltage magnitude and angle before the attack. By minimizing the sum of squared differences between the variable states before and after the attack as the objective function, we can optimize the attack vector in designing the AC FDI attack problem, as represented in Equations \eqref{eq:obj}-\eqref{eq:reactive_overload}.}

\rn{
\begin{small}
\begin{align}
&\min\quad \textstyle\sum_{\small{{m}\in \mathcal{B_A}}}
 (\tilde{V}_{m}- V_{m,fix})^2+ (\tilde{\theta}_{m}-\theta_{m,fix})^2 \label{eq:obj} \\
&\nonumber \text{subject to} \quad \left(\forall i\in\mathcal{B_A}, \forall   \left(l,m\right) \in\mathcal{L_A}\right)\\
&\!\!\!\!\!\!\! g_{sh,i}\, \tilde{V}_i^2+\sum_{\substack{(l,m)\in \mathcal{L},\\\text{s.t.} \hspace{3pt} l=i}} \!\tilde{P}_{lm}+\!\!\sum_{\substack{(l,m)\in \mathcal{L},\\\text{s.t.} 
\label{eq:active_injection}\hspace{3pt} m=i}} \!\!\tilde{P}_{ml}= P_{m,G}-P_{m,D}, \\
&\!\!\!\!\!\!\! -b_{sh,i}\, \tilde{V}_i^2+\!\!\!\!\!\!\sum_{\substack{(l,m)\in \mathcal{L},\\ \text{s.t.} \hspace{3pt} l=i}} \!\!\tilde{Q}_{lm}+\!\!\!\sum_{\substack{(l,m)\in \mathcal{L},\\ \text{s.t.} \hspace{3pt} m=i}} \!\!\!\!\tilde{Q}_{ml}=Q_{m,G}-Q_{m,D},\\
\label{eq:reactive_injection}
&\nonumber\!\!\!\!\!\!\! \tilde{P}_{lm} \!=\! g_{lm} \tilde{V}_l^2\! -\! g_{lm} \tilde{V}_l V_m\cos\left(\tilde{\theta}_{l}-\theta_{m}\right)\!\\
&\qquad\qquad -\! b_{lm} \tilde{V}_l V_m\sin\left(\tilde{\theta}_{l}-\theta_{m}\right),\\
\label{eq:qik1}
&\!\!\!\!\!\!\! \nonumber \tilde{Q}_{lm} = -\left(b_{lm}+b_{c,lm}/2\right) \tilde{V}_l^2 + b_{lm} \tilde{V}_l V_m\cos\left(\tilde{\theta}_{l}-\theta_{m}\right)\\ &\qquad\qquad  - g_{lm} \tilde{V}_l V_m\sin\left(\tilde{\theta}_{l}-\theta_{m}\right),\\
\label{eq:pki1}
&\nonumber \!\!\!\!\!\!\!\tilde{P}_{ml}\! =\! g_{lm} V_m^2\! -\! g_{lm} \tilde{V}_l V_m\cos\left(\tilde{\theta}_{l}-\theta_{m}\right)\!\\
&\qquad\qquad +\! b_{lm} \tilde{V}_l V_m\sin\left(\tilde{\theta}_{l}-\theta_{m}\right),\\
\label{eq:qki1}
&\!\!\!\!\!\!\!\nonumber \tilde{Q}_{ml} = -\left(b_{lm}+b_{c,lm}/2\right) V_m^2 + b_{lm} \tilde{V}_l V_m\cos\left(\tilde{\theta}_{l}-\theta_{m}\right)\\ &\qquad\qquad  + g_{lm} \tilde{V}_l V_m\sin\left(\tilde{\theta}_{l}-\theta_{m}\right)\\
\label{eq:active_overload}
&\!\!\!\!\!\!\! \tilde{P}_{ml} = W_{lm}*{P}_{lm}^{P
F},\\
\label{eq:reactive_overload}
&\!\!\!\!\!\!\! \tilde{Q}_{ml} = W_{lm}*{Q}_{lm}^{P
F}.
\end{align}
\end{small}
}

\rn{In these equations $\tilde{P}_{ml}$, $\tilde{Q}_{ml}$, $\tilde{V}_{l}(\tilde{V}_{m})$ and $\tilde{\theta}_{l}(\tilde{\theta}_{m})$ are respectively active and reactive power flow between busses $m$ and $l$ $(m,l)\in\mathcal{L}$ and voltage values of bussess $m$ and $l$, after the attack. Similarly, $P_{ml}$, $Q_{ml}$, $V_{l}(V_{m})$ and $\theta_{l}(\theta_{m})$ are the same quantities above in before the attack condition. Also, $P_{m,G}$, $Q_{m,G}$, $P_{m,D}$, and $Q_{m,D}$ are active and reactive power generation and active and reactive power demanded in the bus $m$ respectively. It is noticeable that in the Equations~\eqref{eq:active_injection} and~\eqref{eq:reactive_injection} when we write these equations for zero injection busses, they satisfy assumption $1$ in the section~\ref{subsec: AC FDIA}. Also, Equations \ref{eq:active_overload} and \ref{eq:reactive_overload}, are additional constraints to overloading a specified line in the attack zone with a predefined coefficient $W$. By considering this constraint we can design an optimal AC FDI attack for a specific aim like overloading a certain line by a predefined coefficient. ${P}_{lm}^{P
F}$ and ${Q}_{lm}^{P
F}$ are the active and reactive power flow before the attack.}



\rn{After solving this optimization problem, we can calculate the power injection values of non-zero injection busses within the attack zone as follows: }

\begin{subequations}
\begin{small}
\begin{align}
\label{eq:non_zero_injection}
&\tilde{P}_{m}=P_{m}+\sum_{(m,l)\in \mathcal{L}_A}(\tilde{P}_{m,l}-P_{m,l}),\\
&\tilde{Q}_{m}=Q_{m}+\sum_{(m.l)\in \mathcal{L_A}}(\tilde{Q}_{m,l}-Q_{m,l}).
\end{align}
\end{small}
\end{subequations}

\rn{\noindent Here, $\tilde{P_{m}}$, $\tilde{Q_{m}}$, $P_{m}$, and $Q_{m}$ represent the active and reactive power injections at bus $m$ after and before the attack, respectively. Notably, in Equations~\eqref{eq:non_zero_injection}, $\tilde{P}{m,l}$ and $\tilde{Q}{m,l} \in \mathcal{L}_A$ should be considered only for the lines within the attack zone that are connected to bus $m$. Furthermore, when utilizing PMU, the current flow measurement values in the lines should be adjusted based on the corresponding values after the attack. Since the current flows in the lines are related to voltage variables, they can be calculated using Equations~\eqref{eq:from_current_flow} and~\eqref{eq:to_current_flow}, as follows:}

\begin{subequations}
\begin{small}
\begin{align}
\label{eq:from_current_flow}
 &\tilde{I_{f}}= Y_{f}*\tilde{V}\\ 
 \label{eq:to_current_flow}
 &\tilde{I_{t}}= Y_{t}*\tilde{V}
\end{align}
\end{small}
\end{subequations}

\rn{\noindent Here, $\tilde{I_{f}}$ and $\tilde{I_{t}}$ represent the current flow in the lines, from and to buses respectively, while $Y_{f}$ and $Y_{t}$ are the line admittance matrices that account for the admittance of lines from and to the buses, respectively.
After calculating all of these values, we can define the attack vector $a$ based on 
~\eqref{eq:vector a} as follows:}

\begin{small}
\begin{align}
\label{eq:attack_vector_final}
 &\nonumber a= [\tilde{P_{ml}}- P_{ml}, 
     \tilde{Q_{ml}}- Q_{ml},
     \tilde{P_{m}}- P_{m},\\
  &~~~~~~\nonumber   \tilde{Q_{m}}- Q_{m},
     \tilde{V_{m}}\angle \tilde{\theta_{m}}- V_{m}\angle\theta_{m},\\
  &~~~~~~   \tilde{I_{ml}}\angle \tilde{\theta_{ml}}- I_{ml}\angle\theta_{ml},
     \tilde{\theta_{m}}-\theta_{m}]^T
\end{align}
\end{small}

\begin{figure*}
    \centering
    \hspace{-0.8cm}
\captionsetup{justification=centering}
\includegraphics[scale=0.55,trim= 21cm 0.4cm 0.7cm .4cm,clip]{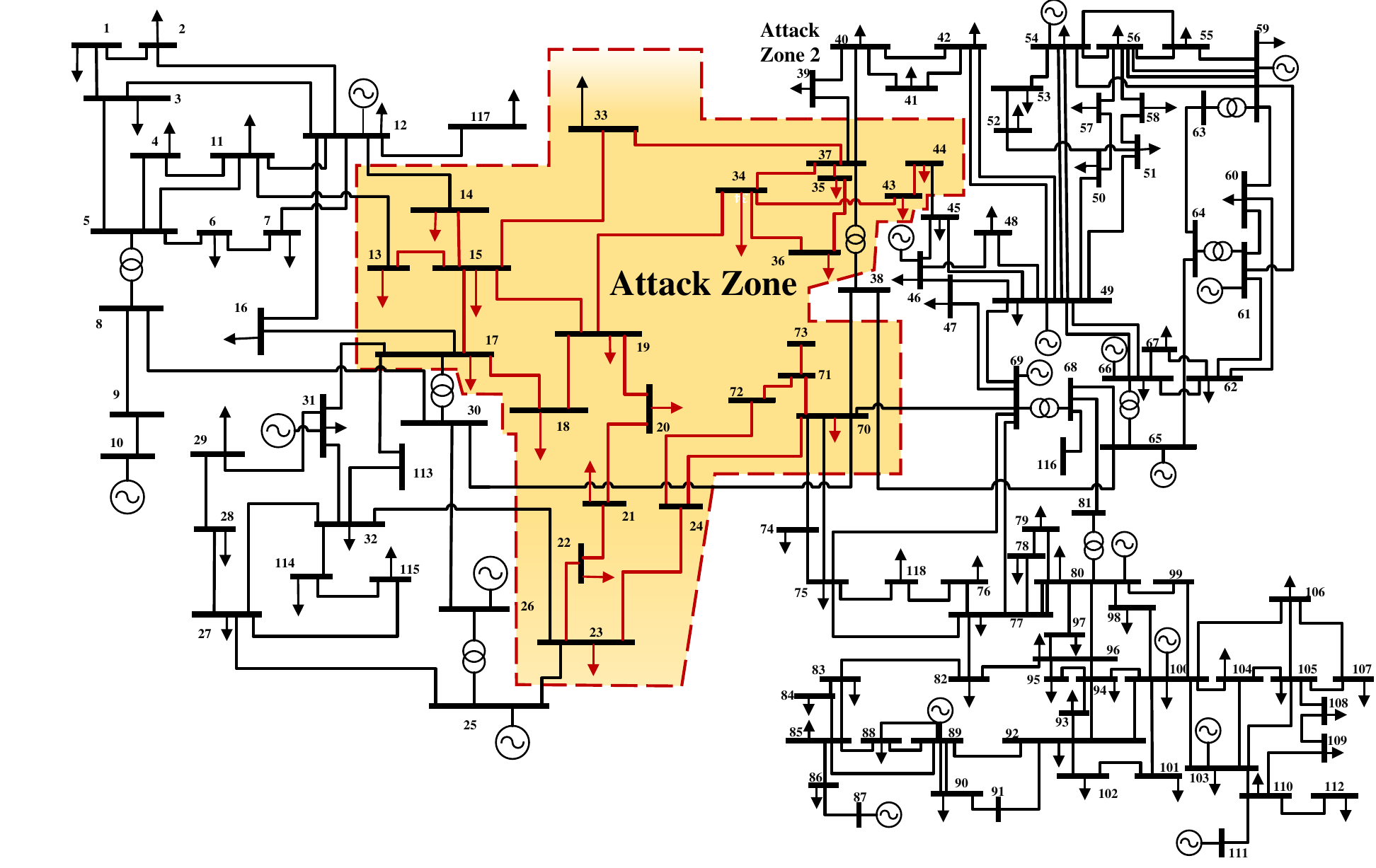}
	 \caption{\rn{One-line diagram depicting the IEEE 118-bus test system, with the attack zones indicated.}}
	\label{fig:attack_zone_118}
\end{figure*}

\section{ANALYSIS OF ATTACK SUCCESS RATE}
\label{sec:ANALYSIS}
\rn{In this section, the proposed approach for designing an optimal AC-FDI attack is applied to the IEEE 118-bus test system from the PGLib-OPF v18.08 benchmark library~\cite{pglib} to test its effectiveness in designing AC FDI attacks. We demonstrate the impact of the designed AC FDI attack on the residuals of estimations. Additionally, we compare optimal and non-optimal FDI attacks to show how the proposed approach makes detecting the designed attack more difficult for traditional residue-based BDD approaches.}

\begin{figure}
    \centering
\includegraphics[scale=0.37]{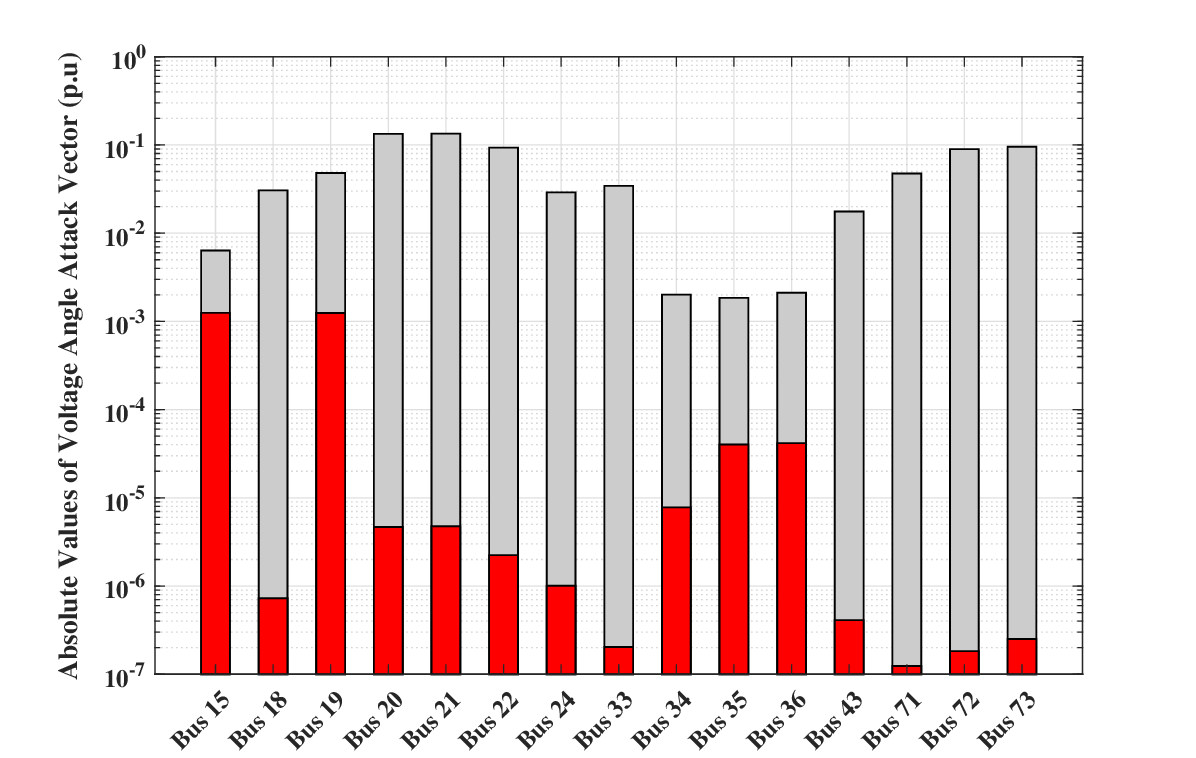}
	\caption{\rn{Attack Vector for voltage angle. The red and gray bars indicate the voltage angles that need to be added to the measurement for optimal and non-optimal FDI attacks, respectively.}}
	\label{fig:voltage_ang118}
\end{figure}

\begin{figure}
    \centering
\includegraphics[scale=0.37]{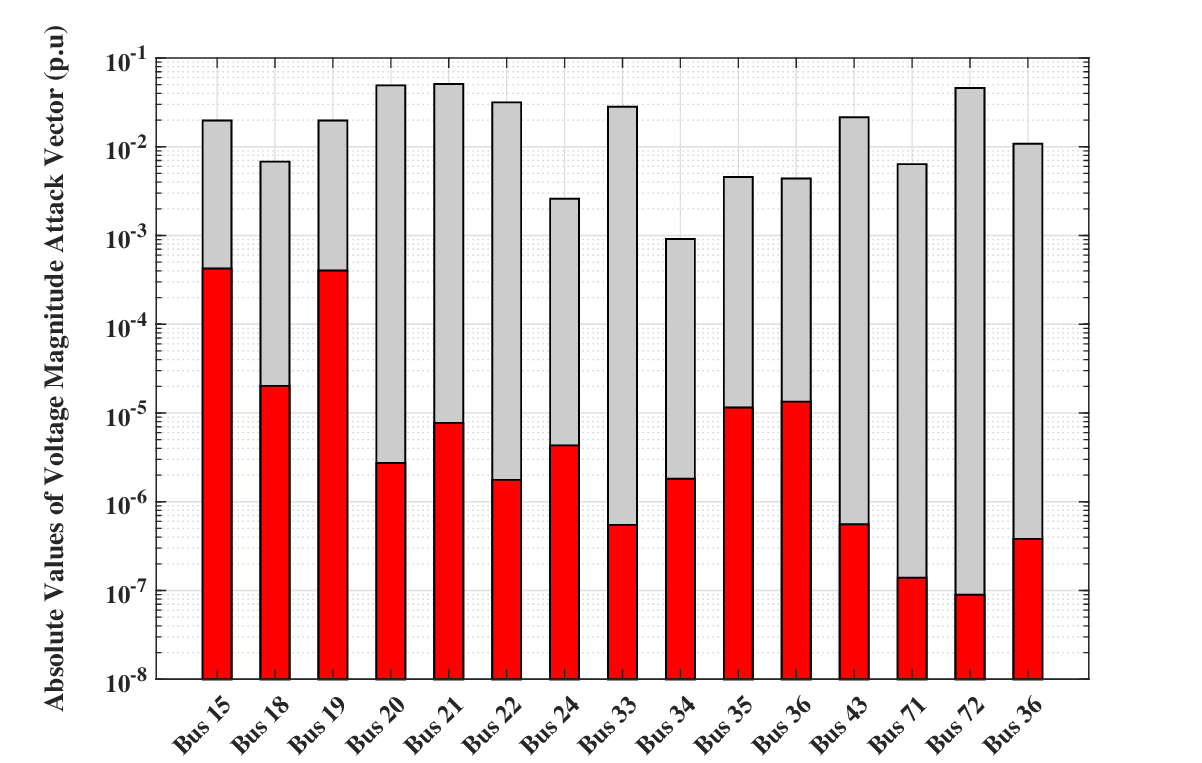}
	\caption{\rn{ Attack Vector for voltage magnitude. The red and gray bars indicate the voltage magnitudes that need to be added to the measurement for optimal and non-optimal FDI attacks, respectively.}}
	\label{fig:voltage_mag118}
\end{figure}

\rn{We considered two scenarios to demonstrate the effects of the proposed attack vector on the residuals of estimations. In the first scenario, we designed the FDI attack without optimizing the attack vector. This was done by substituting the objective function in Equation~\eqref{eq:obj} with a fixed quantity. In this situation, by solving the proposed optimization problem with a fixed objective function, we designed an attack vector ($a$) that satisfied all the constraints required for a successful AC FDI attack. In the second scenario, by considering the objective function in Equation~\eqref{eq:obj} and enforcing the same set of constraints as in the previous scenario, we optimally designed an attack vector ($a$).}

\begin{figure}
    \centering
\includegraphics[scale=0.37]{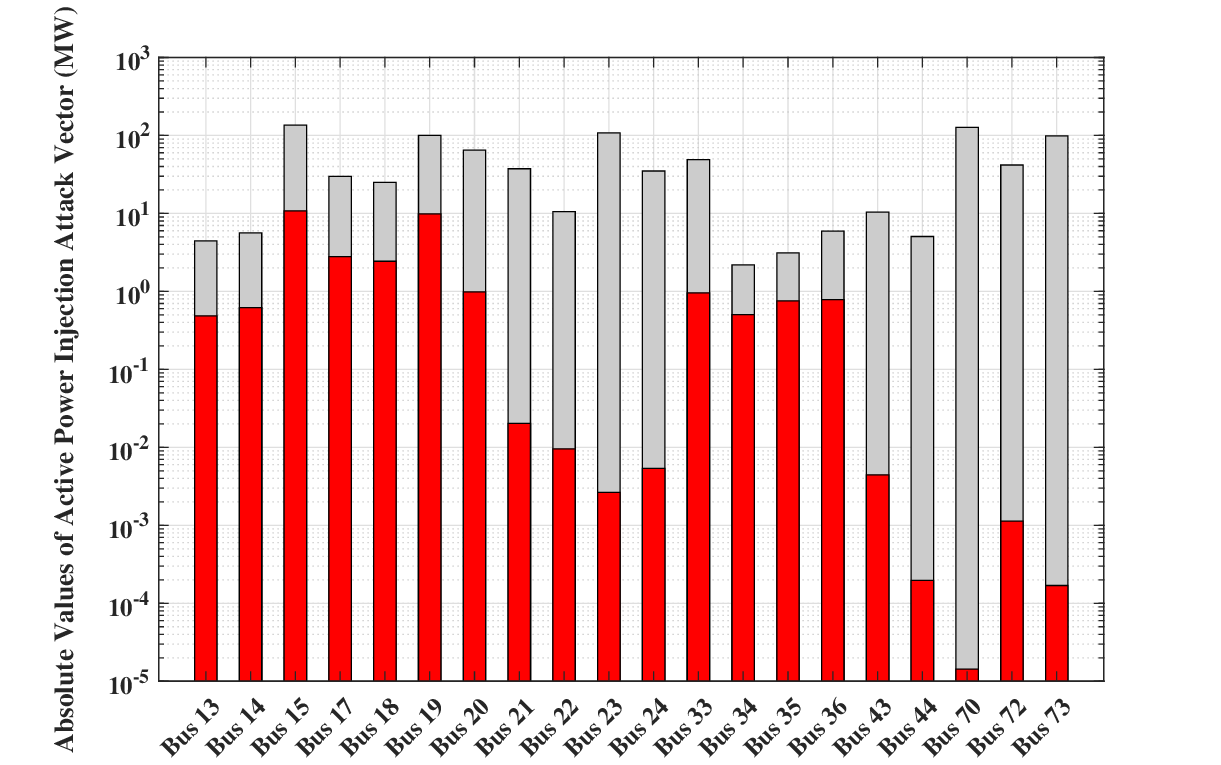}
	\caption{\rn{Attack Vector for active power injection. The red and gray bars indicate the active power injections that need to be added to the measurement for optimal and non-optimal FDI attacks, respectively.}}
	\label{fig:active_injection118}
\end{figure}

\begin{figure}
    \centering
\includegraphics[scale=0.37]{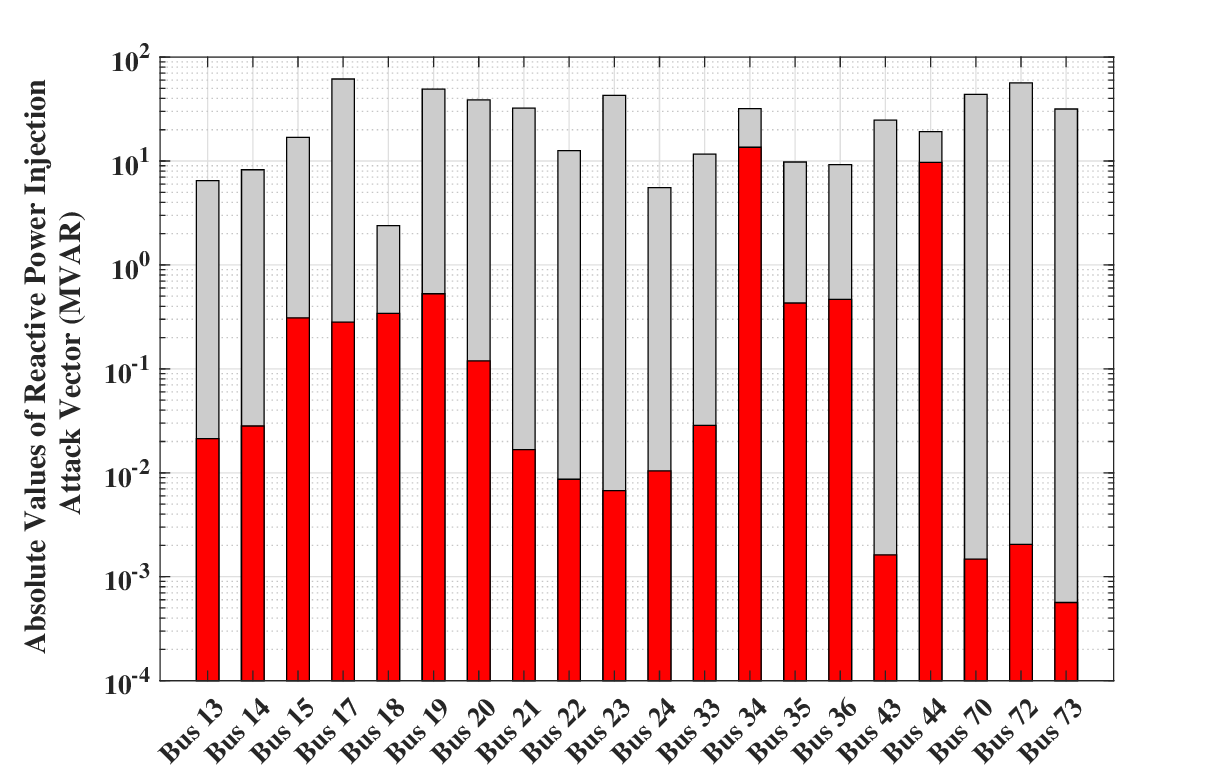}
	\caption{\rn{Attack Vector for reactive power injection. The red and gray bars indicate the reactive power injections that need to be added to the measurement for optimal and non-optimal FDI attacks, respectively.}}
	\label{fig:reactive_injection118}
\end{figure}

\rn{In these two scenarios, as shown in Figure~\ref{fig:attack_zone_118}, we considered 22 buses in the attack zone. Additionally, we aimed to overload line $26$, which connects buses $15$ and $19$, by $1.5$ times its nominal value as a specific goal of this AC-based optimal FDI attack.
Based on this consideration, we have seven boundary buses, including buses $13, 14, 17, 23, 70, 37,$ and $44$. To bind the attack zone, we fixed the voltage magnitude and angle of these boundary buses. Therefore, there are $15$ values for voltage magnitudes and 15 values for voltage angles corresponding to the buses within the attack zone that do not have fixed voltages. These values are shown in Figures~\ref{fig:voltage_ang118} and~\ref{fig:voltage_mag118}.
From these figures, it is evident that the proposed approach for designing an AC FDI attack overloads a line with smaller changes in the magnitude of voltage and angle compared to the non-optimal AC FDI attack approach. Additionally, there are 20 non-zero injection buses in the attack zone. The values of the active and reactive attack vectors for these buses are shown in Figures~\ref{fig:active_injection118} and~\ref{fig:reactive_injection118} for active and reactive power injections, respectively. The changes in the power flow injections of these buses are smaller for the proposed approach compared to the non-optimal approach.}

\begin{figure}
    \centering
\includegraphics[scale=0.37]{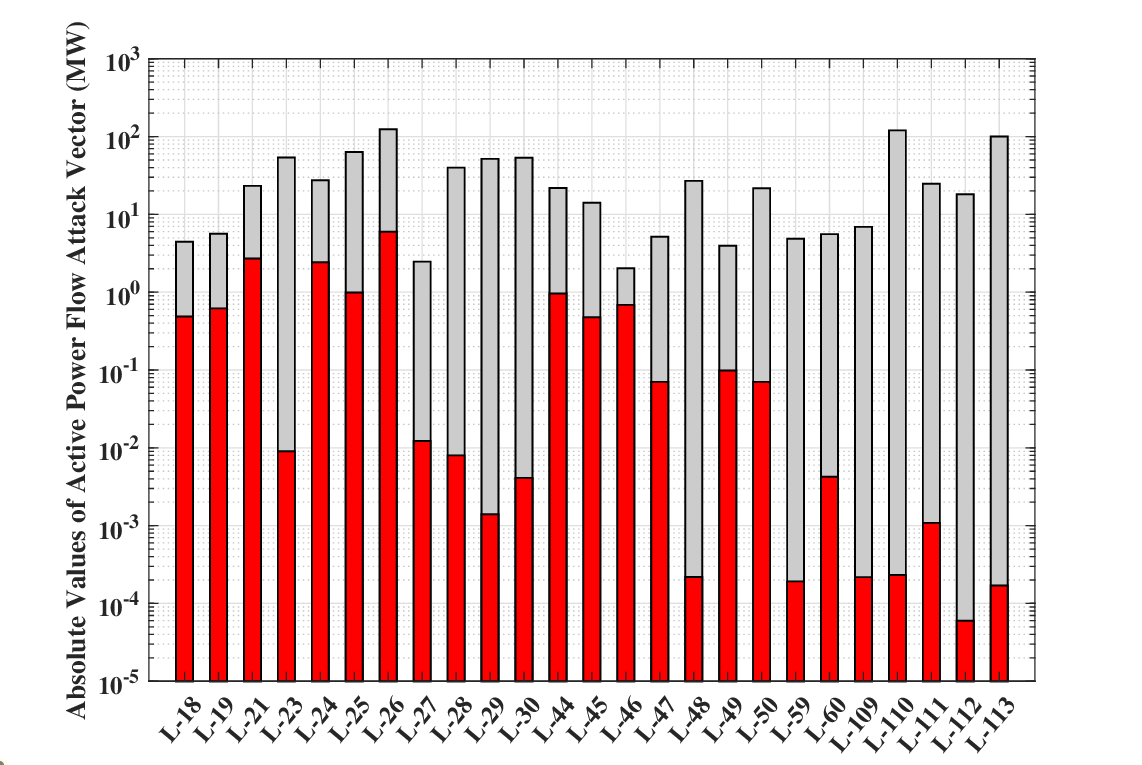}
	\caption{\rn{Attack Vector for active power flow in the lines from the ``sending'' end to the ``receiving'' end. The red and gray bars indicate the active power flows that need to be added to the measurement for optimal and non-optimal FDI attacks, respectively.}}
	\label{fig:active_flow118}
\end{figure}

\begin{figure}
    \centering
\includegraphics[scale=0.37]{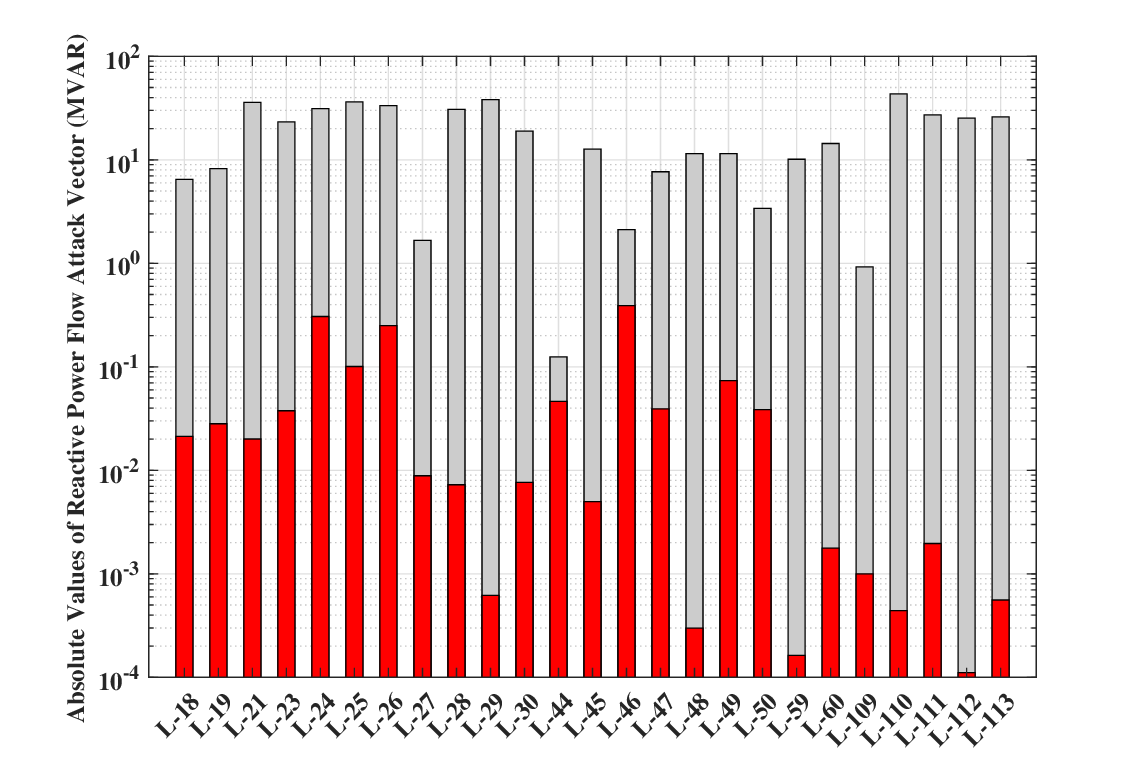}
	\caption{\rn{Attack Vector for reactive power flow in the lines from the ``sending'' end to the ``receiving'' end. The red and gray bars indicate the reactive power flows that need to be added to the measurement for optimal and non-optimal FDI attacks, respectively.}}
	\label{fig:reactive_flow118}
\end{figure}

\begin{figure}
    \centering
\includegraphics[scale=0.37]{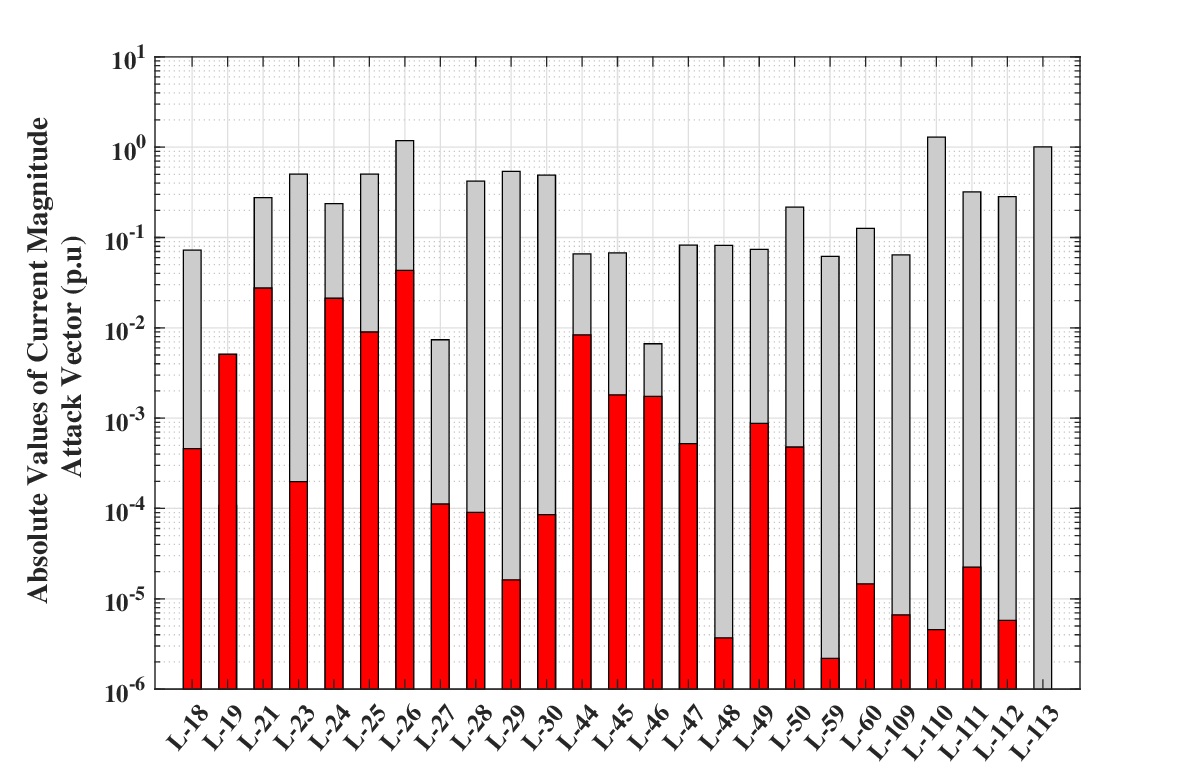}
	\caption{\rn{Attack Vector for current magnitude for different lines. The red and gray bars indicate the current magnitudes that need to be added to the measurement for optimal and non-optimal FDI attacks, respectively.}}
	\label{fig:current_flow118}
\end{figure}

\rn{There are $25$ active power flow attack values for active power flow, $25$ reactive power flow attack values for reactive power flow, and $25$ current magnitude attack values for current flow in $25$ lines in the attack zone of Figure~\ref{fig:attack_zone_118}, which are indicated in red in this Figure. The corresponding values for these quantities are shown in Figures~\ref{fig:active_flow118},~\ref{fig:reactive_flow118}, and~\ref{fig:current_flow118} for both optimal and non-optimal approaches. As can be seen in these figures, the values of the optimal attack vector are smaller than the non-optimal values. In these figures, gray indicates the non-optimal values, and red indicates the optimal values.}

\begin{figure}
    \centering
\includegraphics[scale=0.37]{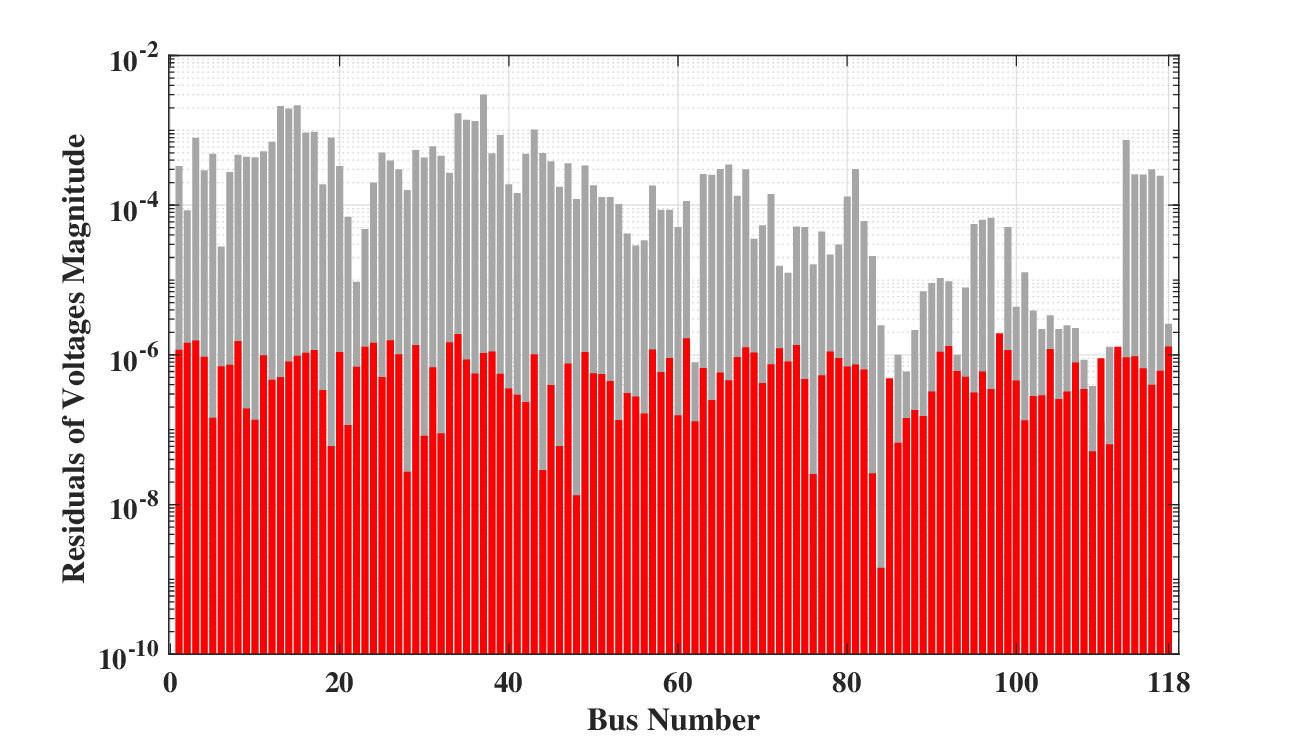}
	\caption{\rn{Residuals of AC state estimator for voltage magnitudes. The red and gray bars indicate the residuals of the AC state estimator for voltage magnitudes under optimal and non-optimal FDI attacks, respectively.}}
	\label{fig:voltage_mag_residuals}
\end{figure}

\begin{figure}
    \centering
\includegraphics[scale=0.37]{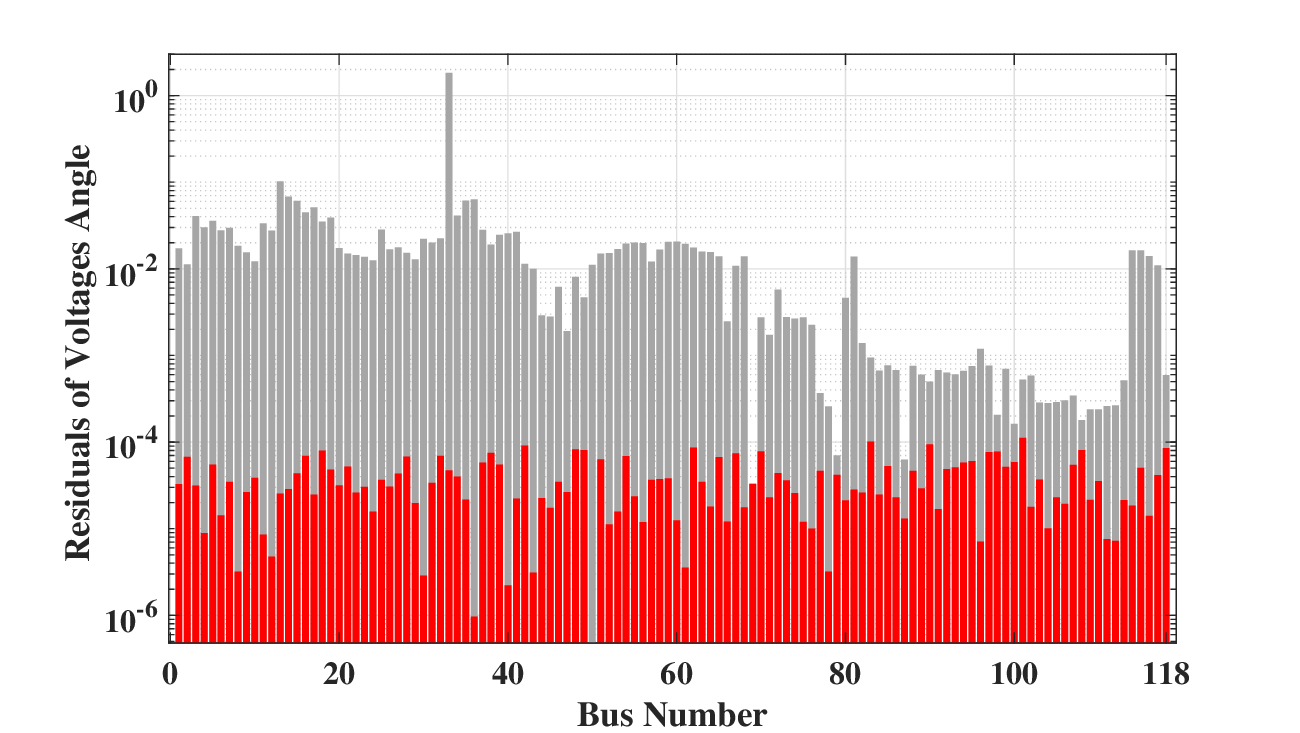}
	\caption{\rn{Residuals of AC state estimator for voltage angles. The red and gray bars indicate the residuals of the AC state estimator for voltage angles under optimal and non-optimal FDI attacks, respectively.}}
	\label{fig:voltage_ang_residuals}
\end{figure}

\begin{figure}
    \centering
\includegraphics[scale=0.37]{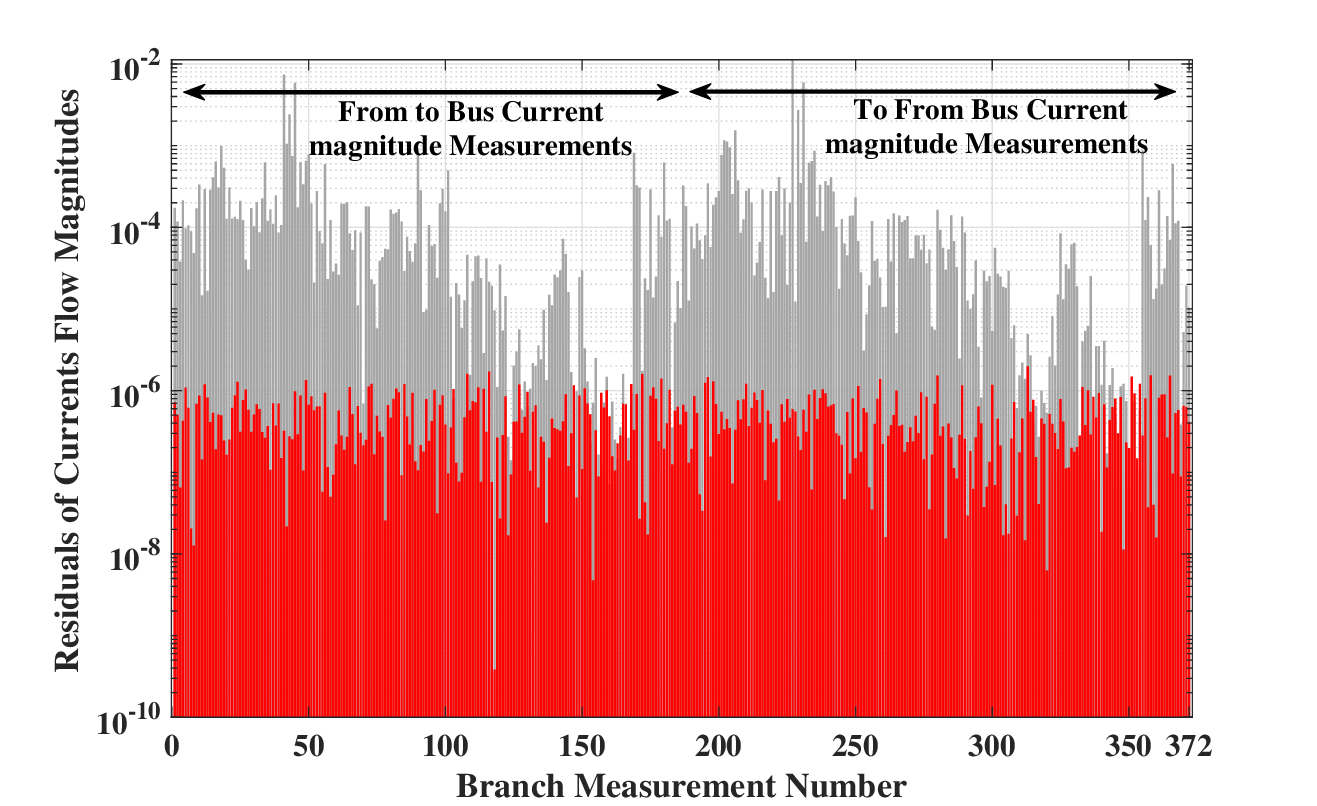}
	\caption{\rn{Residuals of AC state estimator for current magnitudes. The red and gray bars indicate the residuals of the AC state estimator for current magnitudes under optimal and non-optimal FDI attacks, respectively.}}
	\label{fig:current_mag_residuals}
\end{figure}

\rn{To assess the impact of the designed attacks on state estimation, we considered PMU measurements, including voltage magnitude and angle, and current magnitude in the system. By injecting the attack vector into the corresponding measurements and using the MATGRID toolbox~\cite{cosovic2024matgrid} to perform AC state estimation, the residuals for these quantities were calculated and are shown in Figures~\ref{fig:voltage_mag_residuals},~\ref{fig:voltage_ang_residuals}, and~\ref{fig:current_mag_residuals} for voltage magnitude, voltage angle, and current magnitude, respectively. The residuals with optimal attack vectors are significantly smaller than those for non-optimal attack vectors. We demonstrated that although an optimal attack vector could lead to erroneous estimations deviating from the actual system conditions, the small residuals make it challenging to detect using residual-based BDD detection strategies.}

\section{Conclusion}
\label{sec:conclusion}

\rn{This paper proposes a unified approach for designing optimal AC FDI attacks, which can pose challenges for traditional residue-based BDD approaches. We address the system's physics by incorporating power flow equations into the attack design procedure to maintain the undetectability of these proposed attacks by the AC state estimator. Moreover, we define an appropriate objective function to optimize the designed AC FDI attack. Additionally, we impose a constraint on the optimization problem to enforce overloading a line up to a predetermined value. To evaluate the proposed algorithm's ability to design an optimal AC FDI attack, we applied it to the 118-bus test system. The residuals for the AC state estimator with the optimal attack vector are significantly smaller than their counterparts with the non-optimal attack vector. While the proposed AC FDI attack design could lead to erroneous estimations deviating from the actual system conditions, the small residuals in the AC state estimator make it difficult to detect the designed attack by residue-based BDD approaches.}

\bibliographystyle{IEEEtran}
\IEEEtriggeratref{40}
\bibliography{ref}
\end{document}